\documentclass[12pt]{amsproc}
\newtheorem{theorem}{\sc Theorem}[section]
\newtheorem{lemma}[theorem]{\sc Lemma}
\newtheorem{proposition}[theorem]{\sc Proposition}
\newtheorem{corollary}[theorem]{\sc Corollary}

\newtheorem{problem}[theorem]{\sc Problem}

\begin{document}
\title[commutators]
{Multilinear commutators in residually finite groups}
\author{Pavel Shumyatsky}
\address{Department of Mathematics, University of Brasilia,
Brasilia-DF, 70910-900 Brazil}
\email{pavel@mat.unb.br}
\thanks{Supported by CNPq-Brazil}
\begin{abstract} The following result is proved.
Let $w$ be a multilinear commutator and $n$ a positive integer. Suppose
that $G$ is a residually finite group in which every product of at most
896 $w$-values has order dividing $n$. Then the verbal subgroup $w(G)$
is locally finite.
\end{abstract}
\maketitle

\section{Introduction}

According to the solution of the Restricted Burnside
Problem - the celebrated result of Zelmanov \cite{ze1},
\cite{ze2} - every residually finite group of finite exponent is
locally finite. The Lie-theoretic methods used in the solution happened 
to be very effective in treatment of other problems in group theory.
In \cite{shu1999} we used the methods to prove  the following
theorem.

\begin{theorem}\label{1} Let $n$ be a prime-power and
$G$ a residually finite group satisfying the identity
$[x,y]^n\equiv 1$. Then $G'$ is locally finite.
\end{theorem}
Note that in general a periodic residually finite
group need not be locally finite. The corresponding
examples have been constructed in \cite{al,gol,gri,gs,susch}.
Moreover, if the assumption that $G$ is residually finite is dropped
from the hypothesis of Theorem \ref{1}, the derived group need not 
even be periodic. Deryabina and Kozhevnikov showed that for sufficiently
big odd integers $n$ there exist groups $G$ in which all commutators
have order dividing $n$ such that $G'$ has elements of infinite order
\cite{deko}. Independently, this was also proved by Adian \cite{ad}.

 In view of Theorem \ref{1} we raised in \cite{lola} the next problem.

\begin{problem}\label{a} Let $n$ be a positive integer 
and $w$ a word. Assume that $G$ is a residually finite 
group such that any $w$-value in $G$ has order
dividing $n$. Does it follow that the verbal subgroup
$w(G)$ is locally finite?\end{problem}

If $w$ is a word in variables $x_1,\dots,x_t$ we think of
it primarily as a function of $t$ variables defined on any
given group $G$. The corresponding verbal subgroup $w(G)$
is the subgroup of $G$ generated by the values of $w$. The word $w$
is commutator if the sum of the exponents of any variable
involved in $w$ is zero. According to the solution of the Restricted
Burnside Problem the answer to the above question is 
positive if $w(x)=x$. In fact it is easy to see that
the answer is positive whenever $w$ is any 
non-commutator word. Indeed, suppose $w(x_1,\dots,x_t)$
is such a word and that the sum of the exponents of
$x_i$ is $r\neq 0$. Now, given a residually finite
group $G$, substitute the unit for all the variables
except $x_i$ and an arbitrary element $g\in G$ for
$x_i$. We see that $g^r$ is a $w$-value for all
$g\in G$. Hence $G$ satisfies the identity $x^{nr}=1$
and therefore is locally finite by the result of Zelmanov.

Hence, Problem \ref{a} is essentially about commutator
words. In \cite{lola} the problem was solved positively in the case
where $n$ is a prime-power and $w$ a multilinear commutator (outer
commutator word). A word $w$ is a multilinear commutator
if it can be written as a multilinear Lie monomial. Particular examples
of multilinear commutators are the derived words, defined by the equations:
$$\delta_0(x)=x,$$  $$\delta_k(x_1,\dots,x_{2^k})=
[\delta_{k-1}(x_1,\dots,x_{2^{k-1}}),
\delta_{k-1}(x_{2^{k-1}+1}\dots,x_{2^k})],$$
and the lower central words:
$$\gamma_1(x)=x,$$  $$\gamma_{k+1}(x_1,\dots,x_{k+1})=
[\gamma_{k}(x_1,\dots,x_{k}),x_{k+1}].$$

In the case that $n$ is not a prime-power Problem \ref{a} seems to be
very hard. We mention a theorem obtained in \cite{yahahi}.

\begin{theorem}\label{2555} Let $n$ be a positive
integer that is not divisible by $p^2q^2$ for any 
distinct primes $p$ and $q$. Let $G$ be a residually
finite group satisfying the identity
$([x_1,x_2][x_3,x_4])^n\equiv 1$. Then $G'$ is locally finite.
\end{theorem}

Some further progress has been made in \cite{jhone} where the following
theorem was proved.

\begin{theorem}\label{jho}
For any positive integer $n$ there exists $t$ depending only on $n$ 
such that if $w$ is a multilinear commutator and $G$ is a residually
finite group in which every product of $t$ values of $w$ has order
dividing $n$, then $w(G)$ is locally finite.
\end{theorem}

More recently it was shown that in the case that $w=[x,y]$ the number $t$
in the above theorem can be taken to be 68 independently of $n$ \cite{comu68}.
The purpose of the present article is to prove a similar result for arbitrary
multilinear commutators. Thus, we improve Theorem \ref{jho} as follows.

\begin{theorem}\label{main}
Let $n$ be a positive integer and $w$ a multilinear commutator. Let
$G$ be a residually finite group in which every product of 896 $w$-values
has order dividing $n$. Then $w(G)$ is locally finite.
\end{theorem}

The constant 896 in the theorem comes from the famous results of Nikolov
and Segal on commutator width of finite groups. In the course of proving
Theorem \ref{main} we need to consider subgroups of a finite soluble group
that can be generated by 4 $w$-values. In the paper of Segal
\cite{segal} it was shown that every element in the derived group of a
finite soluble $d$-generated group is a product of at most 72$d^2$+46$d$
commutators. A better bound can be obtained working through the proofs
given in \cite{nikolovsegal}. It follows that every element of the derived
group of a finite soluble $d$-generated group is a product of at most
$$\text{min}\{d(6d^2+3d+4),8d(3d+2)\}$$ commutators. In the case that $d=4$
this is 448. We apply this result in the situation where each commutator is
a product of 2 $w$-values. So the constant 896 comes about.

Thus, the theorem of Nikolov and Segal plays an important role in the proof
of Theorem \ref{main}. The proof also relies on the classification of finite
simple groups as well as on the Lie-theoretic techniques that Zelmanov
created in his solution of the Restricted Burnside Problem. We also
use recent result, essentially due to Flavell, Guest and Guralnick, that
an element $a$ of a finite group $G$ belongs to $F_h(G)$  if and only if
every 4 conjugates of $a$ generate a soluble subgroup of Fitting height
at most $h$ \cite{fgg}. It is the necessity to use this result that
accounts for the difference between the constants 896 in Theorem \ref{main}
and 68 in the case of simple commutators $[x,y]$ \cite{comu68}. When
dealing with simple commutators $[x,y]$ it is enough to bound the Fitting
height of some 2-generated subgroups while, as was already mentioned, the case
of general multilinear commutators requires considering 4-generated subgroups.
It should be said that the case of arbitrary multilinear commutators differs
from the case of simple commutators $[x,y]$ in several ways. Probably the most
significant difference occurs when reducing the problem to questions about
finite soluble groups. In the case of commutators $[x,y]$ this was performed
with a relatively short argument as in \cite[Proposition 2.3]{comu68}.
In the present paper the short argument was of no help at all. Instead, we
use a rather intricate Proposition \ref{4} in the next section.

\section{Some useful results on finite groups}

All groups considered in this and the next sections are finite.
We use the expression ``$\{a,b,c,\dots\}$-bounded" to mean ``bounded from
above by some function depending only on $a,b,c,\dots$".  Recall that the
Fitting subgroup $F(G)$ of a group $G$ is the product of all normal
nilpotent subgroups of $G$. The Fitting series of $G$ can be defined by
the rules: $F_{0}(G)=1$, $F_1(G)=F(G)$, $F_{i+1}(G)/F_{i}(G)=F(G/F_{i}(G))$
for $i=1,2,\dots$. If $G$ is a finite soluble group, then the minimal number
$h=h(G)$ such that $F_{h}(G) = G$ is called the Fitting height of $G$.
We will require the following proposition.

\begin{proposition}\label{4conj} Let $G$ be a group and $a\in G$.
Suppose that every subgroup of $G$ that can be generated by four conjugates
of $a$ is soluble with Fitting height at most $h$. Then $a\in F_h(G)$.
\end{proposition}

Essentially, the above proposition is due to Flavell, Guest and Guralnick.
All the tools and arguments that are used in the proof of Proposition
\ref{4conj} can be found in \cite{fgg}. However since the proposition was
not stated in \cite{fgg} explicitly, we outline the proof here.

The fact of crucial importance is that if every four conjugates of an
element $a$ in a group $G$ generate a soluble subgroup, then $a$
belongs to the soluble radical of $G$. This was established independently
in \cite{fgg} and also in Gordeev, Grunewald, Kunyavskii and Plotkin
\cite{ggkp}. The proof uses the classification of finite simple groups.
Therefore it is sufficient to prove Proposition \ref{4conj} under the
additional hypothesis that $G$ is soluble. The key r\^ole in 
the soluble case is played by the following lemma, due to Flavell,
Guest and Guralnick \cite{fgg}.

\begin{lemma}\label{3/4} Let $G$ be a soluble group that possesses an element
$a$ such that $G=\langle a^G\rangle$. Let $k$ be a field. Let V be a
nontrivial irreducible $kG$-module. Then $dim\,C_V(a)\leq\frac{3}{4}\,dim\,V$.
\end{lemma}

The particular case of the above lemma where $a$ has order 3 was established
earlier in Al-Roqi and Flavell \cite{alrflavell}. This enabled the authors
to prove that if $G$ is a soluble group containing an element $a$ of
order 3 such that $G=\langle a^G\rangle$, then there exist four conjugates
of $a$ that generate a subgroup with the same Fitting height as $G$.
Using Lemma \ref{3/4} in place of the particular case of 3-elements in the
Al-Roqi and Flavell arguments we see that the assumption that $a$ is of
order 3 can be dropped and so we obtain the following lemma. 

\begin{lemma}\label{al-roqi} Let $G$ be a soluble group that possesses an
element $a$ such that $G=\langle a^G\rangle$. Then there exist four
conjugates of $a$ that generate a subgroup with the same Fitting height
as $G$.
\end{lemma}

Now the proof of Proposition \ref{4conj} becomes very easy.
\begin{proof} As was mentioned above we can assume that $G$ is soluble.
Set $H_1=\langle a^G\rangle$, $H_2=\langle a^{H_1}\rangle$,
$H_3=\langle a^{H_2}\rangle$, etc. Let $H=\cap_iH_i$. Then $H$ is the
smallest subnormal subgroup of $G$ containing $a$ and it is clear that
$H=\langle a^H\rangle$. It follows that the Fitting height of $H$ is at most
$h$. Since $H$ is subnormal, we conclude that $a\in F_h(G)$.
\end{proof}

We call an element $a$ of $G$ a $\delta_k$-commutator if it is
a value of the word $\delta_k$ in $G$.
A well-known corollary of the Hall-Higman theory \cite{hahi} says
that the Fitting height of a finite soluble group of exponent $n$ is bounded
by a number depending only on $n$. We will denote the number by $h(n)$.

\begin{lemma}\label{heifit} Let $k,n\geq 1$ and $G$ a soluble group in
which every product of 896 $\delta_k$-commutators has order dividing $n$.
Then $h(G)\leq h(n)+k+1$.
\end{lemma}
\begin{proof} Let $a=a_1\in G$ be a $\delta_k$-commutator and $a_2,a_3,a_4$
some conjugates of $a$. Put $H=\langle a_1,a_2,a_3,a_4\rangle$ and
$h=h(n)$. We know from \cite{nikolovsegal} that every element of $H'$ is a
product of 448 commutators of the form $[x,a_i]$ for suitable $x\in G$.
Each commutator of the form $[x,a_i]$ is a product of 2
$\delta_k$-commutators so every element of $H'$ is a product of 896
$\delta_k$-commutators. Hence, $H'$ is of exponent $n$ and so the Fitting
height of $H'$ is at most $h$. It follows that the Fitting height of $H$
is at most $h+1$. We now deduce from Proposition \ref{4conj} that every
$\delta_k$-commutator of $G$ is contained in $F_{h+1}(G)$. Therefore
$G^{(k)}\leq F_{h+1}(G)$ and the lemma follows.
\end{proof}

In what follows $X_k(G)$ denotes the set of all $\delta_k$-commutators
obtained using elements of the group $G$. Note that if $P$ is a subgroup
of $G$, in general $X_k(G)\cap P\neq X_k(P)$.

Let $G$ be a finite group and $k$ a positive integer. We will associate
with $G$ a triple of numerical parameters $n_k(G)=(\lambda,\mu,\nu)$ where
the parameters $\lambda,\mu,\nu$ are defined as follows.

If $G$ is of odd order, we set $\lambda=\mu=\nu=0$. Suppose that $G$ is
of even order and choose a Sylow 2-subgroup $P$ in $G$. 
If the derived length $dl(P)$ of $P$ is at most $k+1$, we define
$\lambda=dl(P)-1$. Put $\mu=2$ if $X_{\lambda}(P)$
contains elements of order greater than two and $\mu=1$ otherwise. We
let $\nu=\mu$ if $X_{\lambda}(P)\not\subseteq Z(P)$ and $\nu=0$ if
$X_{\lambda}(P)\subseteq Z(P)$.

If the derived length of $P$ is at least $k+2$, we define $\lambda=k$. 
Then $\mu$ will denote the number with the property that $2^\mu$ is the
maximum of orders of elements in $X_k(P)$. Finally, we let $2^{\nu}$ be
the maximum of orders of commutators $[a,b]$, where $b\in P$ and $a$ is an
involution in a cyclic subgroup generated by some element from $X_k(P)$.
 
The set of all possible triples $n_k(G)$ is naturally endowed with
the lexicographical order. Following the terminology used by Hall and Higman \cite{hahi}
we call a group $G$ monolithic if it has a unique minimal normal subgroup which is non-abelian simple. In the modern literature such groups very often are called ``almost simple".

\begin{proposition}\label{4} Let $k\geq 1$ and $G$ be a group of
even order such that $G$ has no nontrivial normal soluble
subgroups. Then $G$ possesses a normal subgroup $L$ such that $L$ is
residually monolithic and $n_k(G/L)<n_k(G)$.
\end{proposition}

\begin{proof} Let $M$ be a minimal normal subgroup of $G$. We know that 
$M\cong S_1\times S_2\times\cdots\times S_r $, where $ S_1,S_2,\ldots,S_r$ 
are isomorphic simple groups. The group $G$ acts on $M$ by permuting the
simple factors so we obtain a representation of $G$ by permutations of the
set $\left\{S_1,S_2,\ldots,S_r\right\}$. Let $L_M$ be the kernel of the 
representation. We want to show that $n_k(G/L_M)<n_k(G)$.
Suppose this is not true and $n_k(G/L_M)=n_k(G)$. Let $P$ be a Sylow
2-subgroup of $G$ and assume first that the derived length of $P$ is at
least $k+2$. Suppose further that $\nu(G)\neq0$. Let $q=2^\nu$. Since
$\nu(G/L_M)=\nu(G)$, there exist an involution
$a$ in a cyclic subgroup generated by an element from $X_k(P)$ and $b$
in $P$ such that $[a,b]$ is of order $q$ modulo $L_M$. Then $[a,b]$
permutes regularly some $q$ factors in $\left\{S_1,S_2,\ldots,S_r\right\}$.
Without loss of generality we will assume that $S_1$ is one of those
factors and $S_1,\dots,S_q$ is the corresponding orbit under $[a,b]$. 

Suppose that $a$ takes $S_1$ outside the orbit $S_1,\dots,S_q$.

Let $P_i=P\cap S_i$. Choose a nontrivial element $x\in P_1$ and set
$y=[a,x],c=[a,b]^x$. Then $yc=[a,bx]$ is a commutator of the required form
and therefore $(yc)^q=1$.
Write 
\begin{center}
$1=(yc)^q=yy^{{c}^{-1}}y^{{c}^{-2}}\dots y^c$. 
\end{center}
The element $yy^{{c}^{-1}}y^{{c}^{-2}}\dots y^c$ is a product of
elements of the form $x^{c^j}$, each lying in a different $S_j\in\{S_1,\dots,S_q\}$ and 
elements of the form $x^{-ac^j}$ lying in other simple factors.
Looking at it we conclude that $yy^{{c}^{-1}}y^{{c}^{-2}}\dots y^c\neq 1$.
But that means that the order of $yc$ is divisible by 2$q$, a contradiction.

Therefore for every $i$ the $a$-orbit of $S_i$ is contained in
$S_1,\dots,S_q$. Suppose that
${S_1}^a=S_{i_1},{S_2}^a=S_{i_2},\ldots,{S_q}^a=S_{i_q}$. Again
we look at the expression 
$1=(yc)^q=yy^{{c}^{-1}}y^{{c}^{-2}}\dots y^c$. As above, this is the
product of elements of the form $x^{c^j}$, each lying in a different
 $S_j\in\{S_1,\dots,S_q\}$, and elements of the form $x^{-ac^j}$, each lying in a different
 $S_j\in\{S_1,\dots,S_q\}$ as well. Since both elements $x^{-a},x^{c^{i_1}}$ lie in
$S_{i_1}$ and since $1=(yc)^q$, it follows that $x^{-a}x^{c^{i_1}}=1$ for
every $x\in S_1$. From this and the fact that $c=[a,b]^x$ we deduce that
$ax^{-1}[b,a]^{i_1}$ commutes with $x$. Taking into account that also
$x^{-[b,a]^{i_1}}$ commutes with $x$ (because $x^{-[b,a]^{i_1}}$ belongs
to $S_{q-i_1}$), we conclude that $a[b,a]^{i_1}$ commutes with $x$ for all
$x\in S_1$. Thus, $a[b,a]^{i_1}$ commutes with $S_1$. Similarly we arrive
at the conclusion that $a[b,a]^{i_2}$ commutes with $S_2$,
$a[b,a]^{i_3}$ commutes with $S_3$ and so on. Recall that
${S_1}^{[a,b]}=S_2$, ${S_2}^{[a,b]}=S_3$, $\dots$, etc. Therefore all
the elements
$$[a,b]a[b,a]^{i_2}[b,a],$$ $$[a,b]^2a[a,b]^{i_3}[b,a]^2,$$ $$\vdots$$
$$[a,b]^{q-1}a[a,b]^{i_q}[b,a]^{q-1}$$
commute with $S_1$. Remembering that $a$ is an involution we write
$$[a,b]a[b,a]^{i_2}[b,a]=[a,b]^{i_2+2}a,$$
$$[a,b]^2a[b,a]^{i_3}[b,a]^2=[a,b]^{i_3+4}a,$$
$$\vdots$$
$$[a,b]^{q-1}a[b,a]^{i_q}[b,a]^{q-1}=[a,b]^{i_q+2(q-1)}a.$$
Therefore the elements $[a,b]^{i_1}a,[a,b]^{i_2+2}a,\ldots,[a,b]^{i_q+2(q-1)}a$
commute with $S_1$. Since $[a,b]$ permutes regularly $S_1,S_2,\ldots,S_q$,
we deduce that
$$i_1\equiv i_2+2\equiv i_3+4\equiv\ldots\equiv i_q+2(q-1)\ (\mathrm{mod}\,q).$$
It follows now that $i_1=i_{1+q/2}$ and so $S_1=S_{1+q/2}$,
a contradiction.

Thus, if $\nu(G)\neq 0$, then $\nu(G/L_M)<\nu(G)$. Suppose that $\nu(G)=0$,
that is, every involution $a$ in a cyclic subgroup generated by an element
from $X_k(P)$ is central in $P$. Keeping the above notation 
we remark that since $a$ centralizes $P_i$, it follows that $a$
normalizes $S_i$ for every $i$. Therefore $a\in L_M$ and we conclude
that $\mu(G/L_M)<\mu(G)$.

So $n_k(G/L_M)<n_k(G)$ whenever the derived length of $P$ is at least $k+2$.
We will now assume that the derived length of $P$ is at most $k+1$ and so
$X_\lambda(P)$ is contained in a normal abelian subgroup of $P$. Choose
$d\in X_\lambda(P)$ and $x\in P_i$. Since $[x,d,d]=1$, it is easy to see
that $d$ normalizes $S_i{S_i}^d$. It follows that $d^2\in L_M$. Thus, 
$\mu(G/L_M)=1$. If $\mu(G)=2$, then $n_k(G/L_M)<n_k(G)$ so suppose that
$\mu(G)=1$. If $\nu(G)=0$, it follows that $d\in Z(P)$ and so $d$ centralizes
$P_i$ whence we deduce that $d\in L_M$, in which case
$\lambda(G/L_M)<\lambda(G)$ and we are done. It remains to deal with the
case where $\nu(G)=1$ and $\nu(G/L_M)=1$. In other words, we have to deal
with the case where the elements of $X_\lambda(P)$ are involutions
generating a normal abelian subgroup of $P$. Moreover, there exist
$d\in X_\lambda(P)$ and $b\in P$ such that $[d,b]\not\in L_M$.
It is clear that for some $i$ we have ${S_i}^d\neq {S_i}^{[d,b]}$.
Without loss of generality we assume that ${S_1}^d=S_2$ and
${S_1}^{[d,b]}=S_3$. Now choose $x\in P_1$ and write
$$1=[d,bx]^2=([d,x][d,b]^x)^2=[d,x][d,x]^{[d,b]^x}.$$
This shows that $[d,x]$ commutes with $[d,b]^x$, which easily leads to a
contradiction in view of the assumption that
${S_1}^d=S_2$ and ${S_1}^{[d,b]}=S_3$.

Thus, we have shown that $n_k(G/L_M)<n_k(G)$. Let now $L$ be the intersection
of all the subgroups $L_M$, where $M$ ranges through the minimal normal
subgroups of $G$. It follows that $n_k(G/L)<n_k(G)$ so the proof of the
proposition will be complete once it is shown that $L$ is 
residually monolithic. If $T$ is the product of the minimal normal subgroups
of $G$, it is clear that $T$ is the product of pairwise commuting simple
groups $S_1,S_2,\dots,S_t$ and that $L$ is the intersection of the
normalizers of $S_i$. Since $G$ has no nontrivial normal soluble subgroups,
it follows that $C_G(T)=1$ and therefore any element of $L$ induces a
nontrivial automorphism of some the $S_i$. Let $\rho_i$ be the natural
homomorphism of $L$ into the group of automorphisms of $S_i$. 
It is easy to see that the image of $\rho_i$ is monolithic and that the
intersection of the kernels of all $\rho_i$ is trivial. Hence $L$ is
residually monolithic.
\end{proof}

The next lemma is given without proof as it is precisely Lemma 3.2
from \cite{vari}. The proof is based on Lie-theoretic techniques created
by Zelmanov.

\begin{lemma}\label{2} Let $G$ be a group in which every
${\delta_k}$-commutator is of order dividing $n$. Let $H$
be a nilpotent subgroup of $G$ generated by a set of
$\delta_k$-commutators. Assume that $H$ is in fact $m$-generated
for some $m\geq 1$. Then the order of $H$ is $\{k,m,n\}$-bounded.
\end{lemma}

The lemma that follows partially explains why Proposition \ref{4} is
important for the proof of Theorem \ref{main}.

\begin{lemma}\label{55} There exist  $\{k,n\}$-bounded numbers
$\lambda_0,\mu_0,\nu_0$ with the property that if $G$ is a group
in which every ${\delta_k}$-commutator is of order dividing $n$, then
$n_k(G)\leq(\lambda_0,\mu_0,\nu_0)$.
\end{lemma}
\begin{proof} Suppose that $G$ is a group of even order in which every
${\delta_k}$-commutator is of order dividing $n$ and let $P$ be a
Sylow 2-subgroup of $G$. It suffices to show that there exists a
$\{k,n\}$-bounded number $\nu_0$ such that if $b\in P$ and $a$ is
the involution in a cyclic subgroup generated by some element
$d\in X_k(P)$, then the order of $[a,b]$ is at most $2^{\nu_0}$.
It is clear that $[a,b]\in\langle d,d^b\rangle$. By Lemma \ref{2}
the order of $\langle d,d^b\rangle$ is $\{k,n\}$-bounded so the
result follows.
\end{proof}

\section{Bounding the order of a finite group}

The purpose of this section is to find some sufficient conditions
for a finite group to have bounded order. Recall that all groups
considered in this section are finite.

\begin{lemma}\label{nilpcase} Let $k$, $m$, $n$ be positive
integers and $G$ a nilpotent group in which every
${\delta_k}$-commutator is of order dividing $n$. Assume that $G$ can be
generated by $m$ elements $g_1,g_2,\ldots,g_m$ such that each $g_i$ and
each commutator of the form $[g,x]$, where $g\in\{g_1,g_2,\ldots,g_m\}$
and $x\in G$, have order dividing $n$. Then the order of $G$ is bounded
by a function depending only on $k,m,n$.
\end{lemma}
\begin{proof} Since $G$ is a nilpotent group, it is clear that any
prime divisor of $|G|$ is a divisor of $n$. Hence, it is sufficient to
bound the order of the Sylow $p$-subgroup of $G$ for any prime $p$. We
can pass to the quotient $G/O_{p'}(G)$. Thus, $G$ can be assumed to be
a $p$-group and $n$ a $p$-power.

Suppose first that $G$ is soluble with derived length $j\leq k$. If $G$
is abelian, it is easy to see that $|G|$ is $\{k,m,n\}$-bounded. Arguing
by induction on $j$, we assume that the index $[G:G^{(j-1)}]$ is
$\{k,m,n\}$-bounded.

Suppose that $G^{(j-1)}$ is central. In this case the index $[G:Z(G)]$  is 
$\{k,m,n\}$-bounded and Schur's Theorem \cite[p. 102]{rob} 
guarantees that so is $|G'|$. Since $G$ can be generated by $m$ elements 
of order dividing $n$, it follows that $G$ has $\{k,m,n\}$-bounded order.

Let us see what happens if $G^{(j-1)}$ is not central. Consider the
subgroup
$$A=\langle[g_1,G^{(j-1)}],[g_2,G^{(j-1)}],\ldots,[g_m,G^{(j-1)}]\rangle.$$ 
Clearly, $A$ is normal in $G$. Applying the results of the previous
paragraph to the quotient $G/A$, it follows that $A$ has
$\{k,m,n,\}$-bounded index in $G$. Now Schreier's Theorem says that $A$
can be generated by a $\{k,m,n\}$-bounded number of elements. Since $A$
is abelian and, by the hypothesis, the order of each commutator of the
form $[g,x]$ divides $n$, it follows that $A$ has exponent dividing $n$. 
Hence, $|A|$ and therefore $|G|$ is $\{k,m,n\}$-bounded.

Now consider the general case, that is, we do not assume anymore that
$G$ is soluble with derived length at most $k$. Applying the results of
the previous paragraph to the quotient $G/G^{(k)}$, it follows that
$[G:G^{(k)}]$ is $\{k,m,n\}$-bounded. It remains to show that $|G^{(k)}|$
is $\{k,m,n\}$-bounded. Let $r$ be the minimal number of generators of
$G^{(k)}$. Note that $r$ is $\{k,m,n\}$-bounded. A well-known corollary
of the Burnside Basis Theorem \cite{Huppert1} says that if a $p$-group
is $r$-generated, then any generating set contains a generating set of
precisely $r$ elements. Thus, $G^{(k)}$ can be generated by $r$
$\delta_k$-commutators. By Lemma \ref{2} we conclude
that the order of $G^{(k)}$ is $\{k,r,n\}$-bounded, as required.
\end{proof}

In what follows, we denote by $\pi(G)$ the set of prime divisors of $|G|$.

\begin{lemma}\label{solvcase} Let $k$, $l$, $m$, $n$ be positive
integers and $G$ a group in which every
${\delta_k}$-commutator is of order dividing $n$. Assume that $G$ can be
generated by $m$ elements $g_1,g_2,\ldots,g_m$ such that each $g_i$ and
each commutator of the forms $[g,x]$ and $[g,x,y]$, where
$g\in\{g_1,g_2,\ldots,g_m\}$ and $x,y\in G$, have order dividing $n$.
Assume further that $|G/F(G)|=l$. Then the order of $G$ is bounded by
a function depending only on $k,l,m,n$.
\end{lemma} 
\begin{proof}  Let $F=F(G)$ be the Fitting subgroup of $G$. 
Suppose first that $F$ is central. In this case the index $[G:Z(G)]$
is $\{k,l,m,n\}$-bounded and Schur's Theorem guarantees that so is 
$|G'|$. Since $G$ can be generated by $m$ elements of order dividing $n$, 
it follows that $|G|$ is $\{k,l,m,n\}$-bounded.

If $F$ is not central, consider the subgroup
\begin{center} $N=\langle[g_1,F],[g_2,F],\ldots,[g_m,F]\rangle$.
\end{center}
It is easy to see that $N$ is normal in $G$. Applying the results of the
previous paragraph to the quotient $G/N$, it follows that the index
$[G:N]$ is $\{k,l,m,n\}$-bounded. We will show that $|N|$, and therefore
$|G|$, is $\{k,l,m,n\}$-bounded.

We know that $N$ can be generated by a $\{k,l,m,n\}$-bounded number of
elements. Let $s$ be the minimal number of generators of $N$. Since $N$
is nilpotent, $\pi(N)$ consists of prime divisors of $n$. Thus, it is
sufficient to bound the order of the Sylow $p$-subgroup of $N$ for every
prime $p\in\pi(N)$. Let $P$ be the Sylow $p$-subgroup of $N$ and write
$N=P\times O_{p'}(N)$. If $y_1,y_2,\ldots$ is the list of all elements of
the form $[g_i,y]$, where $1\leq i\leq m$ and $y\in F$, we write
$b_1,b_2,\ldots $ for the corresponding projections of $y_j$ in $P$. Then
$P=\langle b_1,b_2,\ldots\rangle$. Since $P$ is an $s$-generated
$p$-group, the Burnside Basis Theorem shows that $P$ is actually
generated by $s$ elements in the list $b_1,b_2,\ldots$. By the
hypothesis, the order of each of them divides $n$. Each commutator of
the form $[b_i,z]$ also has order dividing $n$. By Lemma \ref{nilpcase}
we conclude that $P$ has $\{k,l,m,n\}$-bounded order.
The proof is complete.
\end{proof}

\begin{proposition}\label{order} Let $k$, $m$, $n$ be positive integers
and $G$ a group in which every product of 896
${\delta_k}$-commutators is of order dividing $n$. Assume that $G$ can be
generated by $m$ elements $g_1,g_2,\ldots,g_m$ such that each $g_i$ and all
commutators of the forms $[g,x]$ and $[g,x,y]$, where
$g\in\{g_1,g_2,\ldots,g_m\}$, $x,y\in G$, have orders dividing $n$.
Then the order of $G$ is bounded by a function depending only on $k,m,n$.
\end{proposition} 
\begin{proof} We use $n_k(G)$ to denote the triple of numerical parameters
as in the previous section. According to Lemma \ref{55} the number of all
triples that can be realized as $n_k(G)$ is $\{k,n\}$-bounded. We therefore
can use induction on $n_k(G)$. If $n_k(G)=(0,0,0)$, then $G$ has odd order.
By the Feit-Thompson Theorem \cite{fetho} $G$ is soluble. By Lemma
\ref{heifit} $h(G)\leq h(n)+k+1$. Arguing by induction on $h(G)$
we can assume that $F(G)$ has $\{k,m,n\}$-bounded index in $G$.
Now the result is immediate from Lemma \ref{solvcase}. Hence, we assume
that $n_k(G)>(0,0,0)$ and there exists a $\{k,m,n\}$-bounded number $N_0$
with the property that if $L$ is a normal subgroup such that
$n_k(G/L)<n_k(G)$, then the index of $L$ is at most $N_0$.

Suppose first that $G$ has no nontrivial normal soluble
subgroups. Proposition \ref{4} tells us that $G$ possesses 
a normal subgroup $L$ such that $L$ is residually monolithic and
$n_k(G/L)<n_k(G)$. It follows that the index of $L$ in $G$ is at most
$N_0$. We deduce that $L$ can be generated by $r$ elements for some 
$\{k,m,n\}$-bounded number $r$.

A result of Jones \cite{Jones} says that any infinite family of finite simple 
groups generates the variety of all groups. It follows that up to isomorphism 
there exist only finitely many monolithic groups in which every 
${\delta_k}$-commutator is of order dividing $n$.
Let $N_1 = N_1(k,n)$ be the maximum of their orders. Then $L$ is
residually of order at most $N_1$. Since $L$ is $r$-generated, the number
of distinct normal subgroups of index at most $N_1$ in $L$ is
$\{r,N_1\}$-bounded \cite[Theorem 7.2.9]{mhall}. Therefore $L$ has 
$\{k,m,n\}$-bounded order. We conclude that $|G|$ is $\{k,m,n\}$-bounded.

Now let us drop the assumption that $G$ has no nontrivial normal
soluble subgroups. Let $S$ be the product of all normal soluble subgroups
of $G$. The above paragraph shows that $G/S$ has $\{k,m,n\}$-bounded
order and we know from Lemma \ref{heifit} that $h(S)\leq h(n)+k+1$.
Again we let $F=F(G)$ be the Fitting subgroup of $G$. Using
induction on the Fitting height of $S$, we assume that $F$ has
$\{k,m,n\}$-bounded index in G, in which case the result is immediate
from Lemma \ref{solvcase}.
\end{proof}

\section{Main results}

We are now ready to prove Theorem \ref{main} in the case where $w$ is a
$\delta_{k}$-commutator.
 
\begin{theorem}\label{deltak} Let $n$ be a positive integer and $G$
a residually finite group in which every product of 896 $\delta_k$-commutators
has order dividing $n$. Then $G^{(k)}$ is locally finite. 
\end{theorem}
\begin{proof}  Let $T$ be any finite subset of $G^{(k)}$. Clearly one can find
finitely many $\delta_{k}$-commutators $h_1,h_2,\ldots,h_m\in G$ such that
$T$ is contained in $H=\langle h_1,h_2,\ldots,h_m\rangle$. It is sufficient to
prove that the subgroup $H$ is finite. The order of each $h_i$ divides $n$.
Moreover, if $h\in\{h_1,h_2,\ldots,h_m\}$ and $x,y\in H$, then the
commutator $[h,x]$ is a product of two $\delta_k$-commutators and the
commutator $[h,x,y]$ is a product of four $\delta_k$-commutators.  So the
order of each of the commutators divides $n$. If $Q$ is any finite quotient
of $G$, by Proposition \ref{order} the order of the image of $H$ in $Q$ is
finite and $\{k,m,n\}$-bounded, so it follows that this order actually does
not depend on $Q$. Since $G$ is residually finite, we conclude that $H$ is
finite, as required.
\end{proof}

Theorem \ref{main} easily follows from Theorem \ref{deltak}. For the
reader's convenience we will reproduce here a couple of lemmas from
\cite{lola}. We say that a multilinear commutator $w$ has weight $k$
if it depends on precisely $k$ independent variables. It is clear
that there are only $k$-boundedly many distinct multilinear commutators
of weight $k$.

\begin{lemma}\label{um} Let $G$ be a group and $w$ a multilinear
commutator of weight $k$. Then every $\delta_k$-commutator in $G$ is
a $w$-value.
\end{lemma}
\begin{proof} The case $k=1$ is quite obvious so
we assume that $k\geq 2$ and use induction on $k$.
Write $w=[w_1,w_2]$, where $w_1$ and $w_2$ are
multilinear commutators of weight $k_1$ and $k_2$
respectively, $k_1+k_2=k$, and the variables 
involved in one of $w_1,w_2$ do not occur in the
other. Let $k_0$ be the maximum of $k_1,k_2$.
By the induction hypothesis any $\delta_{k_0}$-commutator in $G$ is a
$w_1$-value as well as a $w_2$-value. Since $w=[w_1,w_2]$, it follows
that the set of $w$-values contains the set of elements
of the form $[x,y]$, where $x,y$ range 
independently through the set of $\delta_{k_0}$-commutators. Hence any
$\delta_{k_0+1}$-commutator represents a $w$-value. It remains to remark
that $k_0+1\leq k$ so the lemma follows.
\end{proof}
Let $w$ be a multilinear commutator of weight $t$.
In the next lemma we shall require the concept of a
subcommutator of weight $s\leq t$ of $w$. This can be
defined by backward induction on $s$ in the following
way. The only subcommutator of  $w$ of weight $t$  is
$w$ itself. If $s\leq t-1$ a multilinear commutator 
$v$ of weight $s$ is a subcommutator of $w$ if and 
only if there exists a subcommutator $u$ of weight
$>s$ of $w$ and a multilinear commutator $v_1$ such
that either $u=[v,v_1]$ or $u=[v_1,v]$. It is quite
obvious that if $v$ is a subcommutator of $w$ then
$w(G)\leq v(G)$ for any group $G$.
\begin{lemma}\label{dois} Let $w$ be a 
multilinear commutator, $G$ a soluble group in
which all $w$-values have finite order. Then
the verbal subgroup $w(G)$ is locally finite.
\end{lemma}
\begin{proof} Let $G$ be a counter-example whose
derived length is as small as possible, and let
$T$ be the last nontrivial term of the derived
series of $G$. Passing to the quotient over the
subgroup generated by all normal locally finite
subgroups of $G$ we can assume that $G$ has no
nontrivial normal locally finite subgroups. Since
$T$ is abelian, it follows that no $w$-value lies
in $T\setminus\{1\}$. Let $s=s(w,G)$ be the smallest
number such that any subcommutator of weight $s$ of
$w$ has no values in $T\setminus\{1\}$. Obviously,
$s\geq 2$ since $T\neq 1$. We can choose a subcommutator
$v=[v_1,v_2]$ of weight $\geq s$ of $w$ such that
both $v_1$ and $v_2$ are subcommutators of weight
$<s$, at least one of which having nontrivial values
in $T\setminus\{1\}$. Let $H_i$ be the subgroup of $T$
generated by the $v_i$-values contained in $T$;
$i=1,2$. By the choice of $v$ at least one of these
subgroups is nontrivial. Since $v$ has no values in $T\setminus\{1\}$, it
follows that $H_1\leq C_G(v_2(G))$ and $H_2\leq C_G(v_1(G))$. Taking into
account that $w(G)\leq u(G)$ for any subcommutator $u$ of $w$
we conclude that $H_1$ and $H_2$ centralize the verbal subgroup $w(G)$.
Hence both subcommutators $v_1$ and $v_2$ have no nontrivial value in the
image of $T$ in $G/C_G(w(G))$. This shows that $s(w,G/C_G(w(G)))\leq s-1$.
The induction on $s(w,G)$ now shows that $w(G)/Z(w(G))$, the image of $w(G)$
in $G/C_G(w(G))$, is locally finite. Then, by Schur's Theorem, the derived
group of $w(G)$ is locally finite. Because $w(G)$ is
generated by elements of finite order,  $G$  must be locally finite.
\end{proof}

Theorem \ref{main} is now immediate. 
\begin{proof} Indeed, suppose that $G$
satisfies the hypothesis of Theorem \ref{main}. By Lemma \ref{um} there exists 
$k\geq 1$ such that any $\delta_k$-commutator is a $w$-value. Hence any
product of 896 $\delta_k$-commutators in $G$ has order dividing $n$.
Theorem \ref{deltak} now tells us that $G^{(k)}$ is locally finite. It is
straightforward from Lemma \ref{dois} that $w(G)/G^{(k)}$ is likewise locally
finite, as required.
\end{proof}

In \cite{vari} we raised the following problem
that generalizes the Restricted Burnside Problem.

\begin{problem}\label{q1} Let $n\geq 1$ and $w$ a
group-word. Consider the class of all groups $G$
satisfying the identity $w^n\equiv 1$ and having
the verbal subgroup $w(G)$ locally finite. Is that
a variety?
\end{problem}

Recall that variety is a class of groups
defined by equations. More precisely, if $W$ is a
set of words in $x_1,x_2,\dots$, the class of all
groups $G$ such that $W(G)=1$ is called the
variety determined by $W$. By a well-known theorem
of Birkhoff varieties are precisely classes of groups
closed with respect to taking quotients, subgroups and
cartesian products of their members. 

We do not know if Problem \ref{a} and Problem \ref{q1} are equivalent.
It is fairly easy to see that whenever the answer to Problem \ref{q1}
is positive, so is the answer to Problem \ref{a}.
We will show now that for words that are products of multilinear
commutators on independent variables the problems are equivalent indeed.

\begin{proposition}\label{equivale} Let $C$ be a positive integer and
$w$ a multilinear commutator of weight $k$. The following statements are
equivalent.
\begin{enumerate}
\item Every residually finite group $G$ in which all products
of $C$ $w$-values are of order dividing $n$ has $w(G)$ locally finite.
\item Let $G$ be a finite group in which all products of $C$ $w$-values 
are of order dividing $n$. Let $a_1,\dots,a_m$ be $w$-values.
Then the order of $\langle a_1,\dots,a_m\rangle$ is $\{k,m,n\}$-bounded.
\item The class of all groups $G$ in which $w(G)$ is locally finite and
every product of $C$ $w$-values has order dividing $n$ is a variety. 
\end{enumerate}
\end{proposition}
\begin{proof} Suppose first that the first statement is correct but the second
is false. Choose a
family of finite groups $G_1,G_2,\dots,G_i,\dots$ in which all products of
$C$ $w$-values are of order dividing $n$ with the property that for some
$m$ the groups $G_i$ contain $w$-values $a_{i1},\dots,a_{im}$ such that
$$|\langle a_{i1},\dots,a_{im}\rangle|<|\langle a_{j1},\dots,a_{jm}\rangle|$$
whenever $i<j$. Let $G$ be the Cartesian product of the groups $G_i$. It is
clear that $G$ is residually finite. The elements\newline
$b_1=(a_{11},a_{21},\dots,a_{i1},\dots)$\newline
$b_2=(a_{12},a_{22},\dots,a_{i2},\dots)$\newline
$\ldots$\newline
$b_m=(a_{1m},a_{2m},\dots,a_{im},\dots)$\newline
are $w$-values in $G$ but $\langle b_1,\dots,b_m\rangle$ is infinite.
A contradiction. Therefore the first statement implies the second.

Let us now show that the second statement implies the third. We assume that 
the second statement is correct and let $\mathfrak{X}$ denote the class of
all groups $G$ in which $w(G)$ is locally finite and every product of $C$
$w$-values has order dividing $n$. It is easy to see that the class
$\mathfrak{X}$ is closed to taking subgroups and quotients of its members.
Hence, we only need to show that if $D$ is a cartesian product of groups
from $\mathfrak{X}$ then $D\in\mathfrak{X}$. Obviously, every product of
$C$ $w$-values in $D$ has order dividing $n$ so it remains only to prove
that $w(D)$ is locally finite. In view of Lemma \ref{dois} it is
sufficient to show that so is some term of the derived series of $D$.
According to Lemma \ref{um}  every $\delta_k$-commutator in $G$ is
a $w$-value. We wish to show that $D^{(k)}$ is locally finite.
Let $T$ be any finite subset of $D^{(k)}$. Clearly one can find finitely 
many $\delta_{k}$-commutators $h_1,h_2,\ldots,h_m\in D$ such that
$T\leq\langle h_1,h_2,\ldots,h_m\rangle=H$. It is sufficient to
prove that the subgroup $H$ is finite. Since $H$ is generated by finitely many
$\delta_k$-commutators and since every commutator of $\delta_k$-commutators
is again a $\delta_k$-commutator, we deduce that $H^{(k)}$ has finite index
in $H$ and so is generated by finitely many $\delta_k$-commutators, too.
Since the second statement is correct, it follows that the image of $H^{(k)}$
in any finite quotient of $H$ is finite and has bounded order. Thus, it is
sufficient to show that $H$ is residually finite. However this is immediate
from the facts that $H$ is finitely generated and every group $G$ in
$\mathfrak{X}$ has $G^{(k)}$ locally finite.

It remains to show that the third statement implies the first.
Let $\mathfrak{X}$ have the same meaning as in the above paragraph and
assume that $\mathfrak{X}$ is a variety. Let $G$ be a residually finite group
in which all products of $C$ $w$-values are of order dividing $n$.
Then any finite quotient of $G$ belongs to the variety $\mathfrak{X}$. 
However, it is clear that if a group residually belongs to a certain variety, then it 
actually belongs to the variety. Thus, it follows that $w(G)$ is locally finite.
\end{proof}

The next corollaries are now immediate.
\begin{corollary} Let $w$ be a multilinear commutator and $n$ a positive integer.
The class of all groups $G$ in which $w(G)$ is locally finite and
every product of 896 $w$-values has order dividing $n$ is a variety. 
\end{corollary}
\begin{corollary} Let $w$ be a multilinear commutator of weight $k$
and $n$ a positive integer. Let $G$ be a finite group in which every product of
896 $w$-values has order dividing $n$. If $a_1,\dots,a_m\in G$ are $w$-values then
the order of $\langle a_1,\dots,a_m\rangle$ is $\{k,m,n\}$-bounded. 
\end{corollary}

\end{document}